\documentclass{article}
\usepackage[dvips]{graphicx}
\usepackage{amsmath}
\usepackage{amsfonts}
\usepackage[applemac]{inputenc}
\begin{document}
\thispagestyle{empty}

\noindent
\begin{center}{\textbf{On determinacy/indeterminacy of Moment Problems}}\\
\end{center}
\smallskip
\begin{center}\textbf{ERIK ALDÉN}
\end{center}
\bigskip
\noindent
\textbf{0\:\:Introduction}\\

\smallskip
\noindent
This paper treat determinacy of strong moment problems in part I and\\
indeterminacy of strong moment problems in part II.\\
\noindent
This paper is a summary of the following papers:\\

 \noindent      
 [1] \: Aldén, E., Determinancy of Strong Moment Problems.\\
 \indent
\: Department of Mathematics, Umeå University. ISSN 0345-3928\\
\indent
 \: 1987:10, S-901 87 Umeå, Sweden.\\
 
 \noindent     
 [2] \: Aldén, E., On Indeterminancy of Strong Moment Problems.\\
 \indent
 \, Department of Mathematics, Umeå University. ISSN 0345-3928\\
 \indent
  \: 1988:2, S-901 87 Umeå, Sweden.\\
  
\noindent
  [3] \: Aldén, E., Indeterminancy of Strong Moment Problems.\\
 \indent
 \, Department of Mathematics, Umeå University. ISSN 1103-6540\\
 \indent
  \: 1995:7, S-901 87 Umeå, Sweden.\\ 
  
\vskip 0,5cm
\noindent
This paper will treat determinacy/indeterminacy of the strong Stieltjes and Hamburger 
moment problems, part I. Indeterminacy, part II, for certain class of distribution functions. We conclude
by proving a theorem for indeterminacy of the strong problems above for
general distribution functions.\\

\noindent
\textit{Definition 0.1.} \, A function \:$\alpha$\:is called a distribution function, if \:$\alpha(x)$\, is real-valued, bounded and non-decreasing, on some interval \:$I \in R$\,, where \:$R$ \, is the set of real numbers. We also require \:$\alpha$\, to have infinitely many points of increase in \:$I$.\\

\noindent
\textit{Definition 0.2.} \.By a strong moment problem we mean: Given a sequence\\
$\{\mu_n\}_ {n \in Z}$\: of real numbers, find a distribution function \:$\alpha (x)$\, on the interval \:$I$\, such that\\

\begin{equation}
\mu_n = \int_ {I} \,x^n d\alpha(x), n \in Z, \\
\end{equation}
\bigskip

\noindent
where \:$Z$ \, is the set of all integers. The number \:$\mu _n$ \, is called the moment of \:$\alpha $\, of order \:$n$ \,,
$n \in Z$.\\

\noindent
If we take \:$n \in Z_{+}$ \, in \. (1) \, we get the classical moment problem, where \:$Z_{+}$\, is the set of all non-negative integers.\\
The table below explains the various types of moment problems. We may consider an arbitrary bounded interval \:$[a,b]$\, in the Hausdorff case.\\
\begin{center} Hausdorff \hskip 2,5cm Stieltjes \hskip 2,5cm Hamburger
\end{center}
classical \begin{center} $n \in Z_{+}$, $I = [0,1]$ \hskip1cm $n \in Z_{+}$, $I = [0,\infty[$ \hskip1cm $n \in Z_{+}$, $I = ]-\infty,\infty[$
\end{center}
strong \begin{center} $n \in Z$, $I = [0,1]$ \hskip 1cm $n \in Z$, $I = [0,\infty[$ \hskip1cm $n \in Z$, $I = ]-\infty,\infty[$
\end{center}
\smallskip

\noindent
Table 1\\

\noindent
Throughout this paper \:$\alpha$ , \:$\beta$ \, and \:$\sigma$ \, denote distribution functions.\\

\noindent
\textit{Definition 0.3.} A moment problem \:(1) \, is said to be determinate if it has at most one solution \:$\alpha$\, .
Otherwise it is indeterminate.\\

\noindent
\textit{Definition 0.4.} Two distribution functions \:$\alpha_1$\, and \:$\alpha_2$ \, are said to be equal, \:$\alpha_1(x) = \alpha_2(x)$\,, if\\
\begin{equation*} \int_{I} f(x)d\alpha_1(x) = \int_{I} f(x)d\alpha_2(x)
\end{equation*}
\bigskip

\noindent
for all continuous functions \:$f$\, with compact support.\\

\smallskip
\noindent
In the classical Hamburger case,\\

\smallskip
\begin{equation}
\int_{-\infty}^{\infty} \frac{\log\sigma'\,(u)}{1 + u^2} du \ > -\infty,
\end{equation}
\bigskip

\noindent
is a sufficient condition for indeterminacy of \:(1)\, ,see Achiezer [1], p.87.\\
\noindent
We have the corresponding condition for indeterminacy in the classical Stieltjes case,\\

\smallskip
\begin{equation}
\int_{0}^{\infty}\frac{\log\sigma'\,(u^2)}{1 + u^2} du \ > -\infty,
\end{equation}
\bigskip

\noindent
which was proved by the author [3], Theorem 5. We call conditions \:(2)\, and \:(3)\, Krein conditions.\\
\newpage
\noindent
In this paper we will prove that the conditions \:(2)\, and \:(3)\, are sufficient for indeterminacy also of the strong Hamburger and Stieltjes moment problems, provided that the distribution function \:$\sigma(u)$\, has the following
symmetry\\
\noindent
property.\\

\noindent
\textit{Definition 0.5.} A distribution function \:$\sigma$\, on \:$I$\,, where \:$I = ]0,\infty[$\, or \\
$I = ]-\infty,\infty[$\, is called symmetric if\\
\begin{equation}
d\sigma(u) = -d\sigma\Big(\frac{1}{u}\Big)\,, u\in I\setminus \{0\}.
\end{equation}

\bigskip
\noindent
For any symmetric distribution function \:$\sigma$ \, we have \:$\mu_n = \mu_{-n}$\,, $n \in Z$\,, provided that \: $\sigma$ \, has moments of all orders.\\

\noindent
Carleman [7a,7b,7c] gave sufficient conditions for determinacy of the classical Stieltjes and Hamburger moment problems. These conditions are\\
\begin{equation}
\sum_{n=0}^{\infty} \mu_n^{(\frac{-1}{2n})} = \infty
\end{equation}
\noindent
and
\begin{equation}
\sum_{n=0}^{\infty} \mu_{2n}^{(\frac{-1}{2n})} = \infty, 
\end{equation}
\noindent
for the strong moment problems the corresponding results are replaced by\\
\begin{equation}
\sum_{\substack {n\in Z \\ n\neq 0}}\mu_n^{(\frac{-1}{2\vert n\vert})} = \infty 
\end{equation}
\noindent
and
\begin{equation}
\sum_{\substack {n\in Z \\ n\neq 0}} \mu_{2n}^{(\frac{-1}{2\vert n\vert})} = \infty
\end{equation}

\noindent
respectively.\\

\noindent
The conditions \: (5) \, and \: (6) \, are called Carleman conditions. Analogously we call \: (7)\, and \: (8) \, conditions of Carleman type.\\

\noindent
In (1981-84) ,  W.B.  Jones,  O.  Njåstad,  W.J.  Thron and  H.  Waadeland [8,9] have stated and proved necessary and sufficient conditions for the existence for solutions and determinacy. This has been done for both the strong Stieltjes and Hamburger moment problems.\\

\noindent
Two families of strong moment problems are given in section 2. The purpose is to study the sharpness in the sufficient conditions for determinacy.\\
\newpage
\noindent
\textit{Definition 0.6} Let \:$\{\mu_n\}_{-\infty}^{\infty}$ \. be a sequence of non-negative numbers.\\

\noindent
We call the condition\\
\begin{equation*}
\sum_{ \substack {n\in Z \\ n\neq 0}} \mu_n^{(\frac{-1}{2\vert n\vert})} \cdot w_n= \infty, 
\end{equation*}
\bigskip

\noindent
a \underline{weighted Carleman condition} for the strong Stieltjes moment problem \:(1).\\

\noindent
We prove in section 3, that the weighted Carleman condition above with the weight \:$\xi^{\vert n\vert}$ \,is not sufficient for determinacy of the strong Stieltjes moment\\ problem.\\
\noindent
A \textit{family} of strong Stieltjes moment problems\\
\begin{equation*}
\mu_n (d) = \int_{0}^{\infty} x^n d\alpha_d (x) , n \in Z_+
\end{equation*}
\noindent
is a limiting case of the condition of Carleman type \: (7). If\\

\begin{equation*}
\sup_{d} \Bigg \{\sum_ { \substack{n\in Z \\ n\neq 0}} (\mu_n(d))^{(\frac{-1}{2\vert n\vert})}\Bigg \} = \infty,
\end{equation*}
\vskip0,5cm

\noindent\\
Analogously for the conditions \: (2)\,,\: (3)\: and \:(8).\\

\noindent
In Theorem 2.1 we will drop the condition of symmetric distribution function and prove that conditions \.(2)\, and \:(3) \, also are sufficient conditions for indeterminacy in the general strong Hamburger and Stieltjes moment problems respectively.\\

\noindent
The author wants to acknowledge valuable discussions with professor Hans Wallin, fil.dr. Tord Sjödin and fil.dr Per Åhag. Fil.dr Per Åhag has encourage me to write this paper.\\

\newpage
\begin{center}\textbf{Part I}
\end{center}

\bigskip

\noindent
\textbf{1. \:\:Determinacy of the strong Stieltjes moment problem.}
\bigskip

\noindent
The title above may seen a bit obscure, since we shall start to consider\\
 \textit{indeterminacy}. However it vill become clear that indeed we deal with\\
 \noindent
 determinacy.\\

\noindent
This section will deal with a class of strong moment problems a we will prove that this is a limiting case for the conditions of Carleman type as well the condition for indeterminacy, the Krein type condition.\\

\noindent
In [3]  Example 5, p.25, we considered the following family of strong Stieltjes moment problems. Given \:\:
$\{\mu_n(d) \}_{-\infty}^{\infty}$ = $\{\mu_n\}_{-\infty}^{\infty}$ \:\:according to\\
\begin{equation}
\mu_n = \int_{0}^{1} x^n e^{-\big[x^{\frac{1}{(2+d)}} \big]^{-1}} dx +\int_{1}^{\infty} x^n e^{-x^{\frac{1}{(2+d)}} } dx, \: \:n\in Z, \: 0<d\leq 1. \hskip0,8cm  (9)
\end{equation}

\noindent
We found that the moments \:\: $\mu_n$ \, are bounded from above and below by constant multiples of \:\: ${\vert n\vert }^{-(1+\frac{d}{2})}$.
\noindent
The condition of Carleman type \:(7) \, is not fulfilled for the moment problem above but the sums 
\noindent
\begin{equation*}
\sum_{n=0}^{\infty}\mu_n^{(\frac{-1}{2\vert n\vert})} , \\ \sum_{n=-\infty}^{-3}\mu_n^{(\frac{-1}{2\vert n\vert})},\end{equation*}

\noindent
of the positive and negative moments become arbitrarily large when \:\:$d$ \, tends to zero, since for \:\:$d = 0$ \, we have the divergent harmonic series. Hence (9) is determinate for \:\:$d = 0$ .\\
\noindent
When \:\:$d > 0$ \,, we can not make any conclusion about determinacy or\\
indeterminacy.\\
\noindent
From this example we get the idea to construct a family of strong Stieltjes moment problems, by defending a distribution function\:\:$\alpha$ \, with the symmetry property \: (4) \, below.\\
\noindent
If \:\: $\alpha(x)$ \, is a distribution function defined on \:\: $[0,\infty[$ \, then \:\:$-\alpha\big (\frac{1}{x}\big)$ \, also is a distribution function on \:\: $[0,\infty[$.\\
\noindent
\textit{Definition 1.1.}  Suppose \:\: $\alpha(x)$ \, is a distribution function defined on \:\:$ [1,\infty[$ \, and define \:\:$\tilde{\alpha}(x)$ \, according to\\
\begin{equation*}
\tilde{\alpha}(x) = \begin{cases}\-\alpha\big (\frac{1}{x}\big) + m, \: x \in [0,1[, \\ \alpha (x), \:x \in [1,\infty[,\end{cases}
\end{equation*}

\noindent
where \:\:$m$ \, is a real constant such that \:\: $\tilde{\alpha}(x)$ \, becomes non-decreasing. It follows from the construction of \:\: $\tilde\alpha(x)$ \, above that \:\:$\tilde\alpha$ \, is a symmetric distribution function.\\
\noindent
We are now in a position to construct a family of strong Stieltjes moment problems from a family of classical Stieltjes moment problems.\\
\noindent
Let \:\:$\alpha_1$ \, be a distribution function on \:\:$[1,\infty]$ \, with moments of all non-negative orders. Consider the following classical Stieltjes moment problem and corresponding strong moment problem\\
\begin{equation}
\mu_n = \int_0^1 x^n d(-\alpha(\frac{1}{x})) + \int_1^{\infty} x^n d\alpha_1(x) , \, n\in Z_+
\end{equation}

\noindent
and

\begin{equation}
\mu_n = \int_0^1 x^n d(-\alpha(\frac{1}{x})) + \int_1^{\infty} x^n d\alpha_1(x) , \, n\in Z
\end{equation}

\noindent
respectively.\\
\noindent
We now claim that:\\
\noindent
i)\hskip 1cm$\mu_{-n} = \mu_n, \,n\in Z$,\\
ii)\hskip0,9cm if the classical moment problem \:(10) \,is determinate, then the\\
\indent\:\:\:\:\:\:\:\:
corresponding strong Stieltjes moment problem \:(11) \\
\indent \:\:\:\:\:\:\:\: also is determinate.\\
iii)\:\:\:\:\:\:\: \:\:If the classical moment problem \:(10) \, is indeterminate and has\\ 
\indent\:\:\:\:\:\:\:\:\:\,another symmetric distribution function \:\:$\alpha_2$ ,then the strong\\ 
\indent\:\:\:\: \:\:\:\:Stieltjes moment problem \:(11) \, is also indeterminate.\\

\noindent
To prove \:(i) \, we note that\\
\begin{equation*}
\mu_{-n} = \int_0^1 x^{-n} d(-\alpha_1\big(\frac{1}{x}\big)) \,  + \int_1^{\infty} x^{-n} d\alpha_1(x) = \int_0^1x^n d(-\alpha\big(\frac{1}{x}\big)) \,  + \int_1^{\infty} x^n d\alpha_1(x), \,n\in Z ,
\end{equation*}
\hskip1 true cm

\noindent
hence \: (i) \, follows. We omit the proofs of \: (ii) \; and \: (iii) .\\
\newpage
\begin{center}
\textbf{Part II}
\end{center}

\noindent
\textbf{2.\:\:\:Inderminacy of the strong Hamburger moment problem.}\\

\noindent
In this section we will prove that the Krein condition \:(2)\, is a surficient condition for indeterminacy in the strong Hamburger case for a symmetric distribution function.\\

\noindent
\textit{Definition 2.1.} \:\: Let, \:\:$1 \leq p < \infty$ , \, and let \:\:$f$ \, be a real or complex valued and measurable funktion on \:\:$R$ . Let \:\:$\sigma (u)$ , \: $-\infty < u < \infty$ , be a distribution function on \:\:$R$ . We define\\
\begin{equation*}
\lVert{f}\rVert_{p,\sigma} = \Bigg(\int_{-\infty}^{\infty} \lvert f(u)\rvert^p d\sigma(u)\Bigg)^\frac{1}{p}
\end{equation*}

\noindent
and let \:\:${L^p_{\sigma}}$ \, consist of all \:\:$f$ \, for which \:\: $\lVert f \rVert_{p,\sigma} < \infty$.

\noindent
We now claim the following result.\\
\smallskip

\noindent
\textbf{Theorem 2.1.} \:\: Let \:\:$\sigma(u)$ , $-\infty < u < \infty$ ,\, be a symmetric distribution function, and suppose also that \:\:$\sigma(u)$ \, generates finite moments of all orders.\\
\noindent
Let\\
\begin{equation}
\mu_n = \int_{-\infty}^{\infty} u^n d\sigma(u) , n \in Z.
\end{equation}

\noindent
If\\
\begin{equation}
\int_{-\infty}^{\infty} \frac{\log\sigma'\,(u)}{1 + u^2)} du > -\infty,
\end{equation}
\smallskip

\noindent
then the strong Hamburger moment problem \:\:(13) \, is indeterminate and \\$\mu_n = \mu_{-n} , n \in Z$. Here \:\:$\sigma' (u)$ \, is the derivative of the absolutely continuous part of the function \:\:$\sigma(u)$.\\
\smallskip

\noindent
The proof is based on density in \:\:${L^p_{\sigma}}$ \, of the set of all rational functions\\ 
$u^n$ ,$n \in Z$. We can achieve a condition for the density in \:\:${L^p_{\sigma}}$ \, of the set of all rational functions \:\:$u^{-k} , k \in Z_+$ , by considering the linear hull of the functions \:\:$e^{i\alpha\frac{1}{u}} ,\\ \alpha \geq 0$ . This is done in the following lemma, which is the analogue to the lemma in the classical Hamburger case  concerning the linear hull of the functions \:\:$e^{i\alpha u} , \alpha > 0$, see [5] .\\
\newpage
\noindent
\textbf{Lemma 1.} \:\: Suppose that \:\:$\sigma(u) , -\infty < u < \infty$ , is a symmetric distribution function. The linear hull of the functions \:\:$e^{i\alpha\frac{1}{u}} , \alpha \geq 0$ ,\, is dense in \:\:${L^p_{\sigma}} , p \geq 1$ \, if and only if\\
\begin{equation}
\int_{-\infty}^{\infty} \frac{\log\sigma'\,(u)}{1 + u^2)} du = -\infty ,
\end{equation}

\noindent
where\:\:$\sigma'(u)$  \, is the derivative of the absolutely continuous part of the function \:\:$\sigma(u)$.\\
\smallskip

\noindent
\textbf{Proof.} \:\: In view of Krein's theorem \:[1] , p.87, it suffices to prove that the linear hull of \:\:$e^{i\alpha\frac{1}{u}} , \alpha \geq 0$ , is dense in  \:\:${L^p_{\sigma}}$ \, if and only if this holds for the liner hull of \:\:$e^{i\alpha u} , \alpha \geq 0$. First assume that the linear hull of \:\:$e^{i\alpha\frac{1}{u}}$, $\alpha \geq 0$ , is dense in  \:\:${L^p_{\sigma}}$ .
\noindent
Let \:\:$f \in{L^p_{\sigma}} , \, \frac{1}{p} + \frac{1}{q} = 1$ \, be such that\\
\begin{equation*}
\int_{-\infty}^\infty f(u)e^{i\alpha u} d\sigma(u) = 0 ,\, \alpha \geq 0.
\end{equation*}
\smallskip

\noindent
A change of variable gives\\
\begin{equation*}
\int_{-\infty}^\infty  f\big(\frac{1}{u}\big) e^{i\alpha\frac{1}{u}} d\sigma(u) = 0 , \,\alpha \geq 0.
\end{equation*}
\smallskip

\noindent
It now follows that the linear hull of \:\:$e^{i\alpha\frac{1}{u}}, \alpha \geq 0$ , is dense in \:\:${L^p_{\sigma}}$ \,from the Hahn-Banach Theorem and the fact that \:\:$f\big(\frac{1}{u}\big) \in L^p_{\sigma}$. The converse is proved analogously. This proves Lemma 1.\\

\noindent
\textbf{Proof of Theorem 2.1.} \:\:Let \:\:$\sigma$ \, be a distribution function on \:\:$R$ \,satisfying condition \:(14) . Let us first prove that the set of all rational functions\\
\noindent
$u^{-k} , k \in Z_+$ \, is not dense in  \:\:${L^p_{\sigma}} , p \geq 1$ . By Lemma 1 and the Hahn-Banach Theorem there exists a function \:\:$f(u) \in L^p_{\sigma} , f \neq 0$\,,\: where \: $\frac{1}{p} + \frac{1}{q} = 1$ \,, such that\\
\smallskip

\noindent
\begin{equation}
\int_{-\infty}^\infty f(u) e^{i\alpha\frac{1}{u}} d\sigma(u) = 0 , \,\alpha \geq 0.
\end{equation}
\smallskip

\noindent
Now differentiating \:(15) \: $k$ \, times with respect to \:\:$\alpha$ \, and putting \:\:$\alpha = 0$ \, yields\\
\smallskip
\begin{equation*}
\int_{-\infty}^\infty f(u) u^{-k} d\sigma(u) = 0 , k \in Z_+.
\end{equation*}
\smallskip

\noindent
Whence it follows that the set of all rational functions \:\:$u^{-k} , k \in Z_+$ \, is not dense in\:\:$L^p_{\sigma} , p \geq 1$ , by Krein's  theorem, see [1] , p.87. We conclude that the set of all rational functions \:\:$u^n , n \in Z$ , is not dens in \:\:$L^p_{\sigma} , p \geq 1$.\\
\noindent
By the same argument as in the proof in the classical Hamburger case, it follows we have an indeterminate strong Hamburger moment problem,see [1] , p.47 -49.\\
\noindent
This completes the proof of Theorem 2.1.\\

\noindent
\textbf{Example 1.}\:\:Consider the following family of strong Hamburger moment\\
\noindent
problems.\\
\begin{equation}
\mu_n = \mu_{n}(c,d) = \int_{-\infty}^\infty u^n d\sigma(u) = \int_{-\infty}^{\infty} u^n \Bigg(\frac{1}{2}\sqrt{\frac{d}{\pi}}\cdot{\frac{e^{-d(\log \lvert u\rvert)^2}}{{\lvert u \rvert}^c} \bigg)} du,
\end{equation}
\smallskip

\noindent
where \:\:$d \in ]0,\infty[ , c \in R , n \in Z$ , see [4] , p.14.\\

\noindent
(i) \:\: The family \:(16) \,of strong Hamburger moment problems is indeterminate for all \:\,$d$ \: > 0, $ c \in R$ \,, but the distribution function \:\:$\sigma $ \, is symmetric only for \:\:$c = 1$.\\
\begin{equation*}
(ii)\hskip 3cm \mu_n = \begin{cases} \mu_{2k} = e^{\frac{(2k+1-c)^2}{4d}}\, ,\,n = 2k \hskip 0,6cm ,\, k \in Z , \\
\mu_{2k+1} = 0 \hskip 1cm,\, n = 2k + 1\,, k \in Z \end{cases}, \hskip 2cm
\end{equation*}
\smallskip

\noindent
(iii)\:\: the condition \: (13) \, holds for all \:\:$d > 0$ \, and \:\:$c \in R$.\\
\smallskip

\noindent
From this example and Example 2 in the next section, we conclude that there exist indeterminate strong moment problems for which the Krein conditions  \: (2) \,  and \: (3) \, are fulfilled, but the distributen function is not symmetric.\\
\newpage
\noindent
\textbf{3. \:\: Indeterminacy of the strong Stieltjes moment problem.}\\
\smallskip

\noindent
In this section we will prove the corresponding indeterminacy theorem for the strong Stieltjes moment problem for a symmetric distribution function \:\:$\sigma$.\\

\noindent
\textbf{Theorem 3.1.} \:\: Let \:\:$\sigma(u) ,0 \leq u <  \infty$\,, be a symmetric distribution function, and suppose also that \:\:$\sigma(u)$ \, possesses finite moments of all orders.\\
\noindent
Let\\
\begin{equation}
\mu_n = \int_{0}^{\infty} u^n d\sigma (u) , n \in Z.
\end{equation}
\smallskip

\noindent
If\\
\begin{equation}
\int_{0}^{\infty} \frac{\log\sigma'\,(u^2)}{1 + u^2} du > -\infty,
\end{equation}
\smallskip

\noindent
then the strong Stieltjes moment problem \: (17) \, is indeterminate and \\
\noindent
$\mu_n = \mu_{-n} , n \in Z$.\\

\noindent
Theorem 3.1. is proved analogously as in the classical Stieltjes case, with only minor changes, see [5] , Theorem 5. We omit the details.\\

\noindent
\textbf{Example 2.} \:\:( See [4], p-11) Consider the following family of strong Stieltjes moment problems. Let \:\:$d \in [0,\infty[ , c \in R$ \, and\\
\begin{equation}
\mu_n = \mu_{n}(c,d) = \int_{0}^{\infty} u^n d\sigma(u) = \int_{0}^{\infty} u^n \Bigg(\frac{1}{2}\sqrt{\frac{d}{\pi}}\cdot{\frac{e^{-d(\log u)^2}}{{ u }^c}} \bigg) du, n \in Z. \hskip 1,6cm
\end{equation}
\smallskip

\noindent
Then\\
\noindent
(i) \:\: the family of strong stieltjes moment problem \: (19) \, is indeterminate for all \:\:$d > 0 , c \in R$, but the distribution function \:\:$\sigma$ \, is symmetric only for \:\:$c = 1$.\\
\begin{equation*}
(ii)\hskip4,5cm \mu_n = e^\frac{(n + 1 - c)^2}{4d} , n \in Z. \hskip 4cm
\end{equation*}
\smallskip

\noindent
(iii) condition \: (18) \, holds for all \:\:$d > 0$ \, and \:\:$c \in R$.\\
\newpage
\noindent
\textbf{4.\:\:Indeterminacy of the strong Hamburger and Stieltjes moment problems for a distribution function \:\:$\sigma$ \, that is not necessarily of "symmetry type"}\\

\noindent
In this section We will prove that the Krein condition \:(2) \, and \:(3) \, are sufficient conditions for indeterminacy of the strong Hamburger and Stieltjes moment problems. Let us state \textit{Definition 2.1.} once again.\\

\noindent
\textit{Definition 2.1.} \:\: Let, \:\:$1 \leq p < \infty$ , \, and let \:\:$f$ \, be a real or complex valued and measurable funktion on \:\:$R$ . Let \:\:$\sigma (u)$ , \: $-\infty < u < \infty$ , be a distribution function on \:\:$R$ . We define\\
\begin{equation*}
\lVert{f}\rVert_{p,\sigma} = \Bigg(\int_{-\infty}^{\infty} \lvert f(u)\rvert^p d\sigma(u)\Bigg)^\frac{1}{p}
\end{equation*}

\noindent
We will now state our main results.\\

\noindent
\textbf{Theorem 4.1.} \:\: Let \:\:$\sigma(u)$ , $-\infty < u < \infty$ ,\, be a distribution function, and suppose also that \:\:$\sigma(u)$ \, generates finite moments of all orders.\\
\noindent
Let\\
\begin{equation}
\mu_n = \int_{-\infty}^{\infty} u^n d\sigma(u) , n \in Z.
\end{equation}

\noindent
If\\
\begin{equation}
\int_{-\infty}^{\infty} \frac{\log\sigma'\,(u)}{1 + u^2} du > -\infty,
\end{equation}
\smallskip

\noindent
then the strong Hamburger moment problem \:\:(13) \, is indeterminate and. Here \:\:$\sigma' (u)$ \, is the derivative of the absolutely continuous part of the function \:\:$\sigma(u)$.\\

\noindent
The proof is based on non-density in \:\:$L^p_{\sigma}$ \, of the set of all rational\\
functions $\{u^n , n\in Z\}$.
\bigskip

\noindent
\textbf{Proof of Theorem 4.1.}\:\:Suppose that\\

\begin{equation}
 \int_{-\infty}^{\infty} \frac{\log\sigma'\,(u)}{1 + u^2} du > -\infty,
 \end{equation}

\noindent
This is equivalent to the fact that \:\:$e^{i\alpha u}, \alpha \geq 0$ \, is not dens in \:\:$L^p_{\sigma} p\geq 1$ \,, see Achiezer [1] , p. 87. From the Hahn-Banach Theorem there exists a function \:\:$g \in L^{\infty}_{\sigma}$ \, such that\\
\begin{equation}
\int_{-\infty}^{\infty}\ e^{i\alpha u} g(u) d\sigma(u) = 0 , \,\alpha \geq 0.
\end{equation}

\smallskip
\noindent
Differentiating \;\: (23) \:\:\:$n$ \; times with respect to \:\:$\alpha$ \, and putting \:\:$\alpha = 0$ \, yields\\
\begin{equation}
\int_{-\infty}^{\infty}u^n g(u) d\sigma(u)
\end{equation}

\smallskip
\noindent
for \:\:$n \in Z_{+}$. This means that the polynomials \:\:$\{u^n , n\geq 0\}$ \, are non-dense in \:\:$L^p_{\sigma}, p\geq1$. To achieve non-denseness of \:\:$u^n$ \, for negative values of \:\:$n$\\
we proceed as follows. Consider the function \:\:$F_1(u)$ \, defined by\\
\begin{equation*}
\ F_1(u) = \int_{-\infty}^{\infty} \frac{ \ e^{i\alpha u}}{iu} \,g(u) d\sigma(u),
\end{equation*}

\smallskip
\noindent
which is well defined since \:\:$\sigma $ \, has finite moments of all orders. Differentiation of \:\:$F_1(\alpha)$ \, with
respect to \:\:$\alpha$ \, gives \:\:$F_1(\alpha) = 0$\:\,, according to \: (24). From this we conclude that \:\:$F_1(\alpha)$ \, is constant and from the Riemann-Lesbegues lemma the Fourier transform \:\:$F_1\alpha$ \, has the property
\begin{equation*}
\lim_{{n}\rightarrow \infty} F_1(\alpha) = 0
\end{equation*}
\noindent
hence \:\:$F_1(\alpha) = 0$. Hence \:\:$F_1(\alpha ) = 0$.\, Hence  \:\:$F_1(\alpha )$ \, is identically zero and letting \:\:$\alpha$ \, tending to zero gives\\
\begin{equation*}
\int_{-\infty}^{\infty}\frac{1}{u} \,g(u) d\sigma(u) = 0,
\end{equation*}

\noindent
which is \:\:(24) \, for \:\:$n = -1$.\\
\smallskip
\noindent
Now define the function \:\:$F_2(\alpha)$ \, by\\
\begin{equation*}
\ F_2(u) = \int_{-\infty}^{\infty} \frac{ \ e^{i\alpha u}}{(iu)^2}\, g(u) d\sigma(u),
\end{equation*}
\smallskip

\noindent
which is also well defined. Then \:\:$F_2 '\,(\alpha) = F_1(\alpha) \equiv 0$ \, and \:\:$F_2(\alpha) \equiv 0$ \,
by the same argument as above. Letting \:\:$\alpha \rightarrow 0$ \, yields\\
\begin{equation*}
\int_{-\infty}^{\infty}\frac{1}{u^2}\, g(u) d\sigma(u) = 0,
\end{equation*}
\smallskip

\noindent
We continue this process by defining\\

\begin{equation*}
\ F_k(\alpha) = \int_{-\infty}^{\infty} \frac{ \ e^{i\alpha u}}{(iu)^k} \,g(u) d\sigma(u), k \geq 1.
\end{equation*}
\smallskip

\noindent
$F_k(\alpha)$ \, is well defined, since \:\:$\sigma$ \, has finite moments of all orders, and \:\:$F_{k+1}'\,(\alpha) = F_k(\alpha)$\\
\noindent
By an induction argument we can prove that \:\:$F_k(\alpha) \equiv 0$ \, all \:\:$k \geq 1$.\\

\noindent
Now letting \:\:$\alpha$ \, tending to zero we have proved the expression \:(24) \, for \:\:$\alpha$ \, for \:\:$n <  0$ . Since now \:(24) \, holds for all \:\:$n \in Z$ \,, the set of all rational functionns $\{u^n , n\in Z\}$ is not dense in \:\:$L^p_{\sigma}, p\geq 1$ \,,. We are now in a position that we can prove that the strong Hamburger moment problem \:\:(25) \, is indeterminate.\\

\noindent
Let us suppose that \:\:$\lVert{f}\rVert_{{\infty},\sigma} \leq 1$ . Consider the distribution function \:\:$\beta(u)$ \, defined by\\
\begin{equation*}
\beta (u) = \sigma(u)[1 + s\cdot g(u)] , \: u \in R, \:\: s \in [-1,1] .
\end{equation*}
\smallskip

\noindent
It follows \: (24) \,, for \:\:$n \in Z$ \,, that \:\:$\beta (u)$ \, has the same moments as \:\:$\sigma(u)$ \,, and since \:\:$g \not\equiv 0$ \,, \:\:$\beta(u)$ \, and \:\:$\sigma(u)$ \, are different distribution functions possessng the same moments \:\:$\mu_n , n \in Z$. The proof is finished.\\
\bigskip
\noindent
The corresponding result for the Stieltjes moment problem is as follows.\\
\smallskip

\noindent
\textbf{Theorem 4.2.} \:\: Let \:\:$\sigma(u)$ , $-\infty < u < \infty$ ,\, be a distribution function, and suppose also that \:\:$\sigma(u)$ \, possesses finite moments of all orders.\\
\noindent
Let\\
\begin{equation}
\mu_n = \int_{0}^{\infty} u^n d\sigma(u) , n \in Z.
\end{equation}

\noindent
If\\
\begin{equation}
\int_{0}^{\infty} \frac{\log\sigma'\,(u^2)}{1 + u^2} du > -\infty,
\end{equation}
\smallskip

\noindent
then the strong Stieltjes moment problem \:\:(17) \, is indeterminate and. Here \:\:$\sigma' (u)$ \, is the derivative of the absolutely continuous part of the function \:\:$\sigma(u)$.\\
\bigskip

\noindent
Theorem 4.2 is proved as in the strong Hamburger case with only minor changes. See also the proof in the classical Stieltjes case, in \:[3] , p. 10-15.\\
\bigskip

\noindent
\textbf{Proof of Theorem 4.2.} \:\:From [3] , p. 12-13, we have that \:(26) implies that \:\:$\{e^{i\alpha u} , \alpha \geq 0\}$ \, is not dense in \:\:$L_{\sigma}^{\infty}([0,\infty[)$ \, since \:\:$L_{\sigma}^{2} \subset L_{\sigma}^{1}$ . From the Hahn-Banach Theorem there exists a function \:\:$g \in  L_{\sigma}^{\infty}([0,\infty[)$ \, such that\\
\smallskip
\begin{equation}
\int_{0}^{\infty}\ e^{i\alpha u} \,g(u) \,d\sigma(u), \: \alpha \geq 0 .
\end{equation}\\

\noindent
Now we define the function \:\:$F_k(\alpha)$ \, by\\
\smallskip
\begin{equation*}
\ F_k(\alpha) = \int_{0}^{\infty}\frac{ e^{i\alpha u}}{(iu)^k} \,g(u)\, d\sigma(u), k \in Z_+ .
\end{equation*}
\smallskip

\noindent
which is well defined since all moments are finite. First by differentiation of \:\:(27)
with respect to \:\:$\alpha$  , k \, times \:\:$k \geq 0$ \,, and putting \:\:$\alpha = 0$ \, yields that\\

\begin{equation*}
\int_{0}^{\infty} u^k \,g(u) \,d\sigma(u), \: k  \in Z_+ .
\end{equation*}\\
\smallskip

\noindent
By precisely the same argument as in the strong Hamburger case we get that\\

\begin{equation*}
\int_{0}^{\infty}\frac{1}{u^k} \,g(u) \,d\sigma(u), \: k  \in Z_+ .
\end{equation*}
\smallskip

\noindent
and hence \:\:$\{\frac{1}{u^k} , k \in Z_+\}$ \, is not dense in \:\:$L_{\sigma}^{1}([0,\infty[)$. From this we conclude that the set of all rational functions \:\:$\{u^n , n \in Z\}$ \, is not dense in \:\:$L_{\sigma}^{1}([0,\infty[)$.\\

\noindent
Now we may conclude, by the same argument as in the strong Hamburger case, that the Krein condition \:(26) \, is a sufficient condition for indeterminate of the strong Stieltjes moment problem \:(25) \, and the proof is finished.\\
\smallskip

\noindent
\textbf{Remark 4.1.} \:\:Berg [6] , Theorem 5.1, p. 28. has a sufficient condition\\

\begin{equation}
\noindent
2\int_{0}^{\infty}\frac{u \log \sigma  '(u^2)}{1 + u^2} du =  \int_{0}^{\infty}\frac{\log \sigma '(u)}{1 + u^\frac{3}{2}} du > -\infty,
\end{equation}
\noindent
\smallskip

\noindent
for indeterminacy of the classical Stieltjes moment problem.\\

\noindent
This is proved by a more direct approach that the proof the classical Hamburger case in Achiezer but it not known if the condition of Berg above also is a sufficient condition for indeterminacy of the strong Stieltjes moment problem \:\:(25) . It is easy to see that the Krein condition \:(26) \, implies \.(28) \, but it not known if \:(28) \, is a sufficient for the polynomials (rational functions) \:\:$\{u^n , n \in Z_+(n \in Z)\}$ \, to be non-dense in \:\:$L_{\sigma}^{1}(R)$.\\

\noindent
Our main result is that conditions \:(21) \, and \:(26) \, are in fact sufficient conditions for indeterminacy of the strong Hamburger and Stieltjes moment problems, respectively. This is proved in Theorem 4.1 and Theorem 4.2 (Section 4).\\
\noindent
In section 5 we shall define a closely related moment problem, where we construct a symmetric distribution function \:\:$\tilde{\sigma} (u)$ \, from a distribution function \:\:$\sigma (u)$.
\newpage
\noindent
\textbf{5. \:\:The symmetrisized strong Hamburger moment problem}\\

\noindent
Let us first define the following class of distribution functions.\\
\noindent
\textit{Definition 5.1.}\:\:A distribution function \:\:$\sigma$ \, on \:\:$I$ , where\:\:$I = [0,\infty[$ \, or\\
$I = ]-\infty,\infty[$ , is called symmetric if\\
\begin{equation}
\ d\sigma(u) = d\sigma\Big(\frac{1}{u}\Big) , \, u \in I \setminus \{0\}.
\end{equation}
\smallskip

\noindent
For any symmetric distribution function \:\:$\sigma$ \,we have \:\:$\mu_{-n} = \mu_n \,, n \in Z$ \,, provided \:\:$\sigma$ \, has moments of all orders.\\

\noindent
Now consider the strong Hamburger moment problem\\

\begin{equation}
\mu_n = \int_{-\infty}^{\infty} u^n d\sigma(u) , n \in Z.
\end{equation}
\smallskip

\noindent
Suppose also that \:\:$\sigma(u)$ \, possesses finite moments of all orders (hence \:\:$\sigma$ \, has no mass at the origin). Now define  the symmetrisized strong Hamburger moment problem \:(31) \, according to\\

\begin{equation}
\tilde{\mu}_n = \int_{-\infty}^{\infty} u^n d\tilde{\sigma}(u) , n \in Z,  
\end{equation}
\noindent
where\\
\smallskip
\begin{equation}
\tilde{\mu} = \frac{!}{2} \Bigg( \sigma(u) - \sigma\Big(\frac{1}{u}\Big)\Bigg) + Y(\sigma(\infty) - Y(-\infty)) , \: u \in R \setminus\{0\} ,
\end{equation}
\smallskip

\noindent
where \:\:$Y_0$ \, is the Heaviside function located at the origin. Then \:\:$\tilde{\sigma}$ \, is a symmetric distributetion function, in the sense of definition \:(30) \, with finite moments of all orders, see [3] . p.  9-11. The moments \:\:$\tilde\mu_n$ \, are easily calculated. From definition \:(32) \, we get\\

 \begin{equation*}
 \tilde{\mu}_n = \frac{1}{2}(\mu_n + \mu_{-n}) , \, n \in Z.
 \end{equation*} 
\smallskip

\noindent
If \:\:$\sigma$ \, is already symmetric then \:\:$\tilde{\sigma} = \sigma$ \, and \:\:\ $\tilde{\mu}_{n} = \mu_n$ \, $n \in Z$. From Section 4 we know that \: (31) \, is indeterminate if\\

\begin{equation}
\int_{-\infty}^{\infty} \frac{\log\tilde{\sigma}'\,(u)}{1 + u^2} du > -\infty,
\end{equation}
\smallskip

\noindent
The condition \:(32) \, is equivalent to\\
\newpage
\begin{equation}
\int_{-\infty}^{\infty} \frac{\log[\sigma'\,(u) + \frac{1}{u^2}\sigma'\,(\frac{1}{u}u)]}{1 + u^2} du > -\infty.
\end{equation}
\smallskip

\noindent
\textbf{Remark 5.1.} \:\:Since we immediately get the following estimate\\

\begin{equation}
\int_{-\infty}^{\infty} \frac{\log\tilde{\sigma}'\,(u)}{1 + u^2} \,du\, \leq \int_{-\infty}^{\infty} \frac{\log[\sigma'\,(u) + \frac{1}{u^2}\sigma '\,(\frac{1}{u}u)]}{1 + u^2}\, du,
\end{equation}
\smallskip

\noindent
It is obvious that if the Krein condition is fulilled for the moment problem \:(30) \, and hence \:(30) \, is indeterminate, then the strong moment symmetrisized problem Hamburger moment problem \:(31) \,is also indeterminate.\\
\noindent
Except for the relation stated above, in Remark 5.1, it is unknown if an indeterminate strong Hamburger moment problem generally implies that the strong symmetrisized Hamburger moment problem is indeterminate.\\
\noindent
It is also not known if the converse is true.\\
\newpage
\begin{center}
\textbf{References}
\end{center} 
[1]   \:\:Achiezer, N. I., 1963, The Classical Moment Problem and Some Related Questions \\
\indent
 \:\:in Analysis. Oliver and Boyd, Edinburgh, London.\\
 \vskip 0,8cm
 \noindent
[2]   \:\:Achiezer, N. I., 1945, On a Proposition of Kolmogorov and a Suggestion of Krein.\\
\indent
\:\:Dokl. Akad. Nauk. SSSR, \textbf{50} (in Russian).\\
\vskip 0,8cm
\noindent
[3]   \:\:Aldén, E., 1985, A survey of weak and strong moment problems with  \\
\indent
\:\:generalizations, University of Umeå, No. 2. \\
\vskip 0,8cm
\noindent
[4]   \:\:Aldén, E., 1987, Determinacy of strong moment problems,  \\
\indent
\:\,University of Umeå, No. 10. \\
\vskip 0,8cm
\noindent
[5]   \:\:Aldén, E., 1988, On Indeterminacy of Strong Moment Problems.  \\
\indent
\:\,Department of Mathematics, University of Umeå, No. 2. \\
\vskip 0,8cm
\noindent
[6]   \:\:Berg, C., 1994, Indeterminate Moment Problems and the Theory of Entire \\ 
\indent
\:Functions. K\o benhavs Universitet, Matematisk Institut, Preprint Seriies, No. 27.\\
\vskip 0,8cm
\noindent
[7a] \:Carleman, T., 1960, Edition complète des articles de Torsten Carleman.\\
\indent
\:\:L' institut mathématique Mittag-Leffler avec l' appui du conseil national Suedois\\
\indent
\:\,pour la recherches dans les sciences naturelles.\\
\indent
\vskip 0,8cm
\noindent
[7b] \:Carleman, T., 1922, Sur le Problème des moment. \\
\indent
\:\:Copmples rendus, 174, 1680- 1682.\\
\vskip 0,8cm
\noindent
[7c] \:Carleman, T., 1923, Sur les équations intégrales singulières a noyau réel et \\
\indent
\:\:symétrique. Uppsala Universitets Årsskrift, 228 pp.\\
\indent
\vskip 0,8cm
\newpage
\noindent
[8]   \:\: Jones, W. B., Thron, W. J.,  Waadeland, H., 1980, The strong Stieltjes moment \\
\indent
\:\:problem. Trans. Amer. Math. Soc. 261, 2, 503- 528.\\
\vskip 0,8cm
\noindent
[9]   \:\: Jones, W. B., Thron, W. J.,  Njåstad, O., 1984, Orthogonal Laurent \\
\indent
\:\:Polynomials and the strong Hamburger moment problem. Journal of \\
\indent
\:\:Mathematical Analysis and Journal of Mathematical Analysis and applications.\\
\indent
\:\:V. 98, No. 2.\\
\vskip 0,8cm
\noindent
[10]  \:\:Shohat, J. A., Tamarkin, J. D., 1943, The Problem of Moments. \\ 
\indent
\:\: Amer. Math. Soc. Providence, R. I.\\ 
\vskip 0,8cm
\noindent
[11]   \:\:Sjödin, T., 1986, A note on the Carleman condition for determinacy of moment \\
\indent\
\:\:problems. University of Umeå, No. 1.\\
\vskip 0,8cm
\noindent
[12]   \:\:Stieltjes, T. J., 1982, Recherches sur les fractions continues. \\
\indent\
\:\:Mémoires presentés par divers Savants à l' académie des Sciences de \\
\indent
\:\: l' institut de France, Science et Mathematiques, (2), 32, (2), 196 pp.\\
\vskip 0,8cm
\noindent
[13]   \:\:Stieltjes, T. J., 1894, Recherches sur les fractions continues. \\
\indent\
\:\:Anales de la Faculté des Sciences de Toulouse, (1), 8, T1- 122, (1), 9, A5- A7. \\
\indent
\end{document}